\documentclass[twocolumn]{cinc}
\usepackage{graphicx}
\usepackage{amsmath}
\usepackage{amssymb}
\usepackage{xcolor}
\usepackage{comment}
\begin{document}
\bibliographystyle{cinc}

\title{Assessment of a Threshold Method for Computing Activation Maps from Reconstructed Transmembrane Voltages}

\author {Emma Lagracie$^{1}$, Lisl Weynans$^{1}$, Yves Coudière$^{1}$ \\
\ \\ 
$^1$ Univ. Bordeaux, INRIA, IMB, UMR 5251, IHU Liryc, F-33400 Talence, France}

\maketitle

\begin{abstract}
Electrocardiographic imaging non-invasively reconstructs activation maps of the heart from temporal body surface potential maps by post-processing solutions of an inverse problem. Typically, activation times are detected through the maximal deflection of the temporal or spatial derivative of recovered extracellular or transmembrane potentials. However, this method can introduce artificial lines of block in the map, falsely indicating a pathology. Consequently, several complex algorithms have been developed in an attempt to smooth activation maps while preserving true discontinuities. 

We propose a straightforward method for computing activation maps from recovered transmembrane voltages, wherein activation time is defined as the first time a predefined threshold is crossed. We evaluate this threshold-based method against traditional deflection-based methods using simulated data, following the approach of Shuler et al. \cite{schuler2021reducing}.

Our findings indicate that the threshold method, when combined with classical Tikhonov regularization, produces smooth activation maps even in the presence of true lines of block. Given the variability in performance of deflection-based methods compared to those reported in \cite{schuler2021reducing}, we emphasize the difficulty of establishing a universally effective post-processing method for computing activation maps.
\end{abstract}


\section{Introduction} 
Electrocardiographic imaging (ECGi) provides a non-invasive image of the heart's electrical activity based on the electrical potentials of the body surface. In particular, a purpose of ECGi is to retrieve activation maps of a heart, i.e. the juxtaposition of the activation times of each point on the heart volume or surface. Activation maps are a powerful tool for cardiologists, as they can indicate pathologies such as areas of low conduction or abnormal pacing sites.

They are usually built by solving an inverse problem in terms of electrical potentials or voltages from body surface data at several time steps, then post-processing the reconstructed signals to extract activation times. A regularization term must be added to constrain and stabilize the solutions of the inverse problem because it is ill posed. However, this regularization can significantly modify the behavior of the recontructions (excessive smoothing, oscillations, etc.). From the reconstructed signals, activation times are generally detected as the time of maximum deflection on the temporal or spatial derivative of the extracellular potential $U$.

Authors of \cite{duchateau2019performance,schuler2021reducing} have shown that the combination of the regularization and post-processing techniques may generate artificial lines of block (LOBs) in the computed activation maps. However, these lines may reveal pathology such as areas of low conduction. There is therefore great interest in reducing these artificial block lines, while still being able to detect true discontinuities. In \cite{schuler2021reducing}, the authors showed that only the combination of the temporal and spatial derivative of the transmembrane voltage $V$ with strong temporal smoothing avoids the appearance of artificial LOBs while preserving true ones. 

According to the stereotyped shape of the transmembrane voltage $V$, the activation time can be considered as the time of maximum derivative (in time and space), but also as the first moment when a threshold value is crossed. This value can be chosen as the mid value between the rest phase and the plateau phase, or above. In this article, we extend the work of \cite{schuler2021reducing} by comparing this new method called “threshold method” with the “Defl-ST” spatio-temporal derivation method from \cite{schuler2021reducing} applied on extracellular potential $U$ and transmembrane voltage $V$. 

\section{Method}
We simulate body surface potential map (BSPM) data using classical forward electrophysiology models. Then, using a source model developed in \cite{nous} and a different mesh, we compute the solutions $U$ and $V$ of the inverse problem on the heart surface. Three post-processing methods are applied to the recovered $V$ and $U$ and compared.
\subsection{Forward simulations}
\begin{itemize}
    \item[\textbf{Source model}] The bidomain model coupled with the Mitchell-Shaeffer ionic model \cite{mitchell2003two} was used to compute the transmembrane voltage $V$, the extracellular and extracardiac potential $U$, and the BSPM data. The intracellular, extracellular and extracardiac conductivity tensors were chosen isotropic for a matter of simplicity. 
    \item[\textbf{Mesh and numerical method}] We used a tetrahedral epi-endocardial mesh (about 20000 nodes), surrounded by a volume torso mesh (approximately 1500 body surface nodes). The transfer matrices were computed using the finite element method.
    \item[\textbf{Propagation scenarios}] Different activation propagation scenarios, of which three relevant are outlined in this article (Figure \ref{FIGURA1}), were initiated by a temporally localized stimulation current over a small area of the heart. In Case 1, the propagation of the activation sequence of the heart is initiated in the septum, while in Cases 2 and 3, the earliest activation sites are on the ventricles. Moreover, in Cases 2 and 3, we inserted low conductivity regions, provoking lines of block in the reference activation maps. 
\end{itemize}

\subsection{Inverse reconstruction}

\begin{itemize}
    \item[\textbf{Source model}] The source model is based on the "averaged model" developed in a previous work \cite{nous}, that couples the transmembrane voltage and extracellular potential on the epicardium with the extracardiac potential in the torso volume. We briefly explain the model construction. The electrical signals across a small thickness are averaged, and thus reduced to surface signals. Here, we take into account a 5 mm deep layer of myocardium from the pericardial surface, and consider the averaged transmembrane voltage and extracellular potential across that thin layer, very close to the true epicardial signal. This model allows to incorporate naturally the heart-torso fluxes as a ratio of the in-depth variation of the extracardiac potential accross the layer. Incorporating those fluxes is very important to reproduce the behaviour of a non-isolated heart. Moreover, using this model allows to recover the transmembrane voltage $V$ simultaneously to the extracellular potential $U$ on the epicardium. The blood fluxes are neglected. 
    
    In the case of isotropic conductivities, the model writes
     \begin{equation*}
\label{ITF}
    \begin{aligned}
        & (\sigma_i + \sigma_e) \Delta_S U + \frac{\sigma_T}{h} \frac{\partial U_T}{\partial n }|_{\partial H} + \sigma_{i} \Delta_S V = 0 &  \text{ on }\partial H,\\
        & \sigma_T \Delta U_T = 0 & \text{ in } T,
    \end{aligned}
\end{equation*}with the boundary conditions
\begin{equation*}
\label{BC}
    \begin{aligned}
        &\sigma_T \partial_n U_T  = 0 & \text{on } \Gamma_T, \\
        & \sigma_{e} \frac{U_T- U}{\alpha h} = \sigma_T \partial_n U_T & \text{on } \partial H.
    \end{aligned}
\end{equation*} The domains $\partial H$, $T$ and $\Gamma_T$ denote respectively the epicardium, the torso volume and the body surface. The conductivities $\sigma$ are indexed by their domains: $i$ accounting for intracellular, $e$ for extracellular and $T$ for torso. We denote by $h$ the depth of the myocardial layer, and $\alpha$ a fixed coefficient between $0$ and $1$.
    \item[\textbf{Mesh and numerical method}] We used a pericardial surface mesh for the heart (about 17000 nodes), coupled with a volume mesh for the torso (about 1200 body surface nodes). Again, we used the finite element method (FEM) for computing the transfer matrix. 
    \item[\textbf{Data}] The BSPM data from the forward simulations were interpolated onto the inverse mesh, and a 4\% Gaussian noise was added to produce the input data $Z_T$ of the inverse reconstruction. We used a different realization of the noise at each time step.
    \item[\textbf{Formulation of the inverse problem and regularization}]\, 
    
    For each time step associated to a BSPM data $Z_T$, the solutions $\Bar{U}$ and $\Bar{V}$ of the inverse problem satisfy 
    $$ (\Bar{U}, \Bar{V}) = \text{arg min} \|A (U, V) - Z_T \| ^2 + \varepsilon \| \sigma_i \Delta V \| ^2,
    $$ where $A$ is the FEM transfer matrix associated with the (inverse) source model, and $\| . \|$ is the FEM $L^2$ norm. The term $\| \sigma_i \Delta V \| ^2$ corresponds to the second order Tikhonov regularization on $V$. The regularization parameter $\varepsilon$ was empirically chosen at $10^{-2}$.
\end{itemize}

\subsection{Post-processing techniques}
\subsubsection{Computing the activation maps}
\begin{itemize}
    \item[\textbf{Post-processing of the transmembrane voltage}] At each \\
    time step, the transmembrane voltage $V$ can only be recovered up to a constant. So, to obtain a coherent time course, this constant is adjusted as in \cite{schuler2019delay}. Before the time of maximal deflection $t_d$ (beginning of the activation process), the minimum of $V$ is fixed at $0$, while after $t_d$ ($\forall~  t>t_d$), the maximum of $V(t)$ is fixed at the value $\text{Max} \, V(t_d)$. 
    \item[\textbf{Temporal smoothing}] The temporal signals $U$, $V$, $\nabla U$ and $\nabla V$ were smoothed with a Gaussian filter of order 0 and standard deviation of 10 ms. 
    \item[\textbf{Spatio-temporal derivative}] The spatio-temporal derivative is defined as the product of the temporal derivative and the norm of the surface gradient of a signal. To compute the surface gradient of $U$ and $V$, we used properties of the Lagrange $P^1$ finite elements, that are affine in each triangle of the mesh. From the Tikhonov regularization, the gradients are naturally smooth. The temporal derivatives were computed using a centered Euler scheme on the temporally smoothed $U$ and $V$. The "Defl-ST" activation times were defined as the time of maximal positive spatio-temporal derivative for $V$, and the time of maximal negative spatio-temporal derivative for $U$. 
    \item[\textbf{Other deflection-based methods}] Only the Defl-ST method is presented in the paper, but other deflection-based methods were implemented and compared. The Defl-ST method produced the more accurate deflection-based maps, avoiding some artificial LOBs.
    \item[\textbf{Threshold method}] In the threshold method, activation times are defined as the first time at which $V$ crosses a certain threshold, located at the mid value between the plateau and rest values of the action potential. However, to limit the impact of the oscillations and amplitude variations due to the Tikhonov regularization, we computed several activation maps with thresholds uniformly distributed between $\frac{2}{3} V_\text{max} - \frac{1}{6} V_\text{max} $ and $\frac{2}{3} V_\text{max} + \frac{1}{6} V_\text{max}$. Then we averaged these maps for obtaining a final activation map.
\end{itemize}

\subsubsection{Error measurements}
We define three types of error to measure the accuracy of an activation map. They are the $L^2$ relative error (\textbf{L2}), the Pearson's Correlation Coefficient (\textbf{CC}), and the Slowness Coefficient (\textbf{SC}). The relative $L^2$ error corresponds to a root mean square error, except that the punctual errors are weighted by the local mesh size:
$$\text{\textbf{L2err}} = \sqrt{\frac{\left (\text{AT} - \text{AT}_\textbf{Ref}  \right )^T \mathbb{M} \left (\text{AT} - \text{AT}_\textbf{Ref}  \right ) }{ \text{AT}_\textbf{Ref}^T \mathbb{M} \text{AT}_\textbf{Ref}}}.$$ The matrix $\mathbb{M}$ is the finite element mass matrix over the heart surface. The Pearson's correlation coefficient for a data serie $X$ writes $$\textbf{CC} = \frac{\sum ^n _{i=1}(X[i] - \overline{X})(X[i] - \overline{X}_\textbf{Ref})}{\sqrt{\sum ^n _{i=1}(X[i] - \overline{X})^2} \sqrt{\sum ^n _{i=1}(X_\textbf{Ref}[i] - \overline{X}_\textbf{Ref})^2}},
$$ where $\overline{X}$ stands for the average value of $X$. For the \textbf{CC} error of the activation map, $X$ is replaced by the $\text{AT}$. The Slowness Coefficient was introduced in \cite{schuler2021reducing} and is defined by the correlation coefficient between the components of the surface gradient of the activation map of reference, and those of the computed activation map. This indicator has the advantage of being more sensible to lines of block than the Correlation Coefficient.

\section{Results}
\begin{figure*}[htbp]
\centering
\includegraphics[width=11cm]{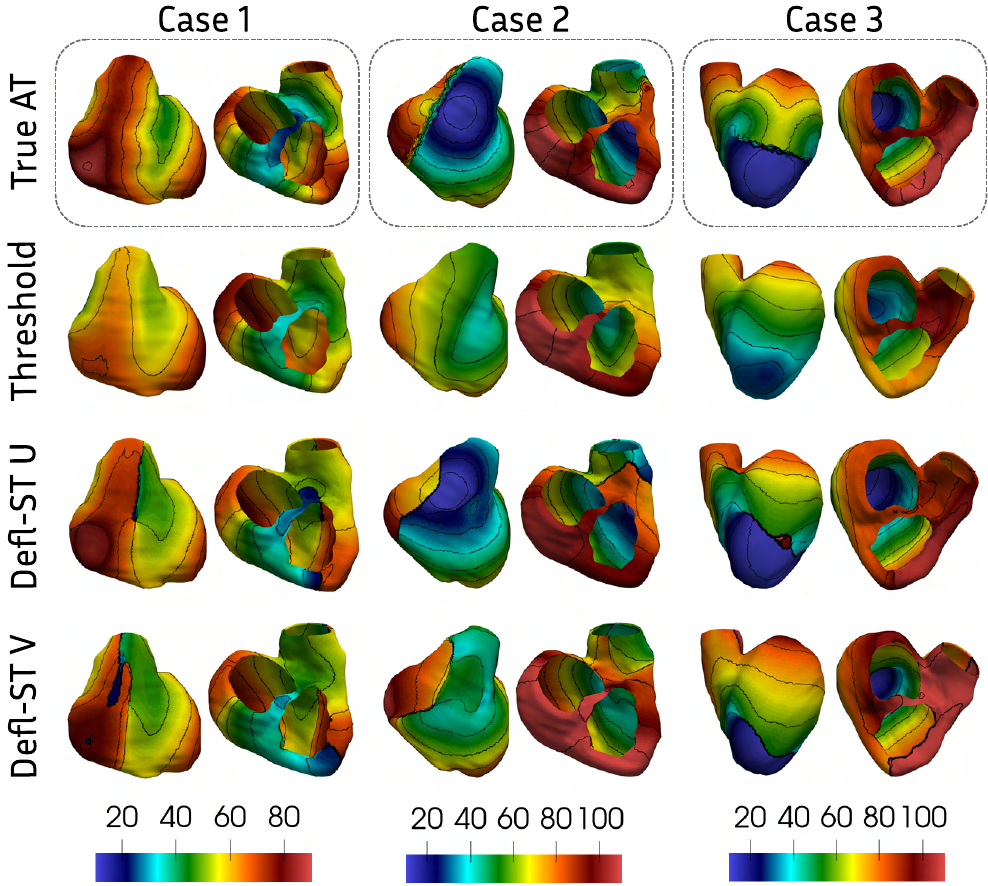}
\caption{Activation maps (ms) computed with the three methods (Threshold, Defl-ST U, Defl-ST V) for three propagation scenarios (Cases 1 - 3).}
\label{FIGURA1}
\end{figure*}

\begin{figure*}[htbp]
\centering
\includegraphics[width=14.7cm]{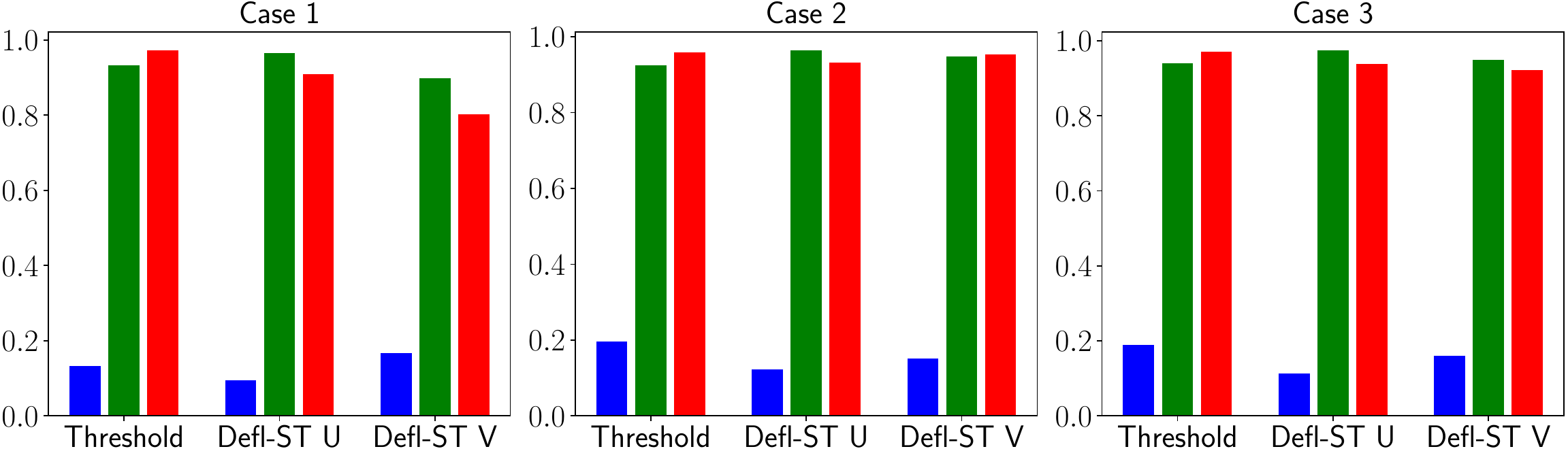}
\caption{Errors on the activation maps computed with the three methods for the scenarios Cases 1 - 3. The blue bar indicates the relative \textbf{L2} error, the green bar the correlation coefficient (\textbf{CC}) and the red bar the slowness coefficient (\textbf{SC}).}
\label{fig:err}
\end{figure*}

Results are presented on Figures \ref{FIGURA1} and \ref{fig:err} through three examples that correspond to the three propagation scenarios. 
We observe that the threshold method provides smoother activation maps than the deflection based methods, and thus avoids the formation of artificial lines of block (Case 1). However, in Cases 2 and 3, the threshold method does not reproduce the real LOBs that are present in the reference AT maps (see Figure \ref{FIGURA1}). Quantitatively, the slowness coefficient \textbf{SC} of the threshold maps is higher than that of the Defl-ST maps in the three cases, but the correlation coefficient remains below or equal to the Defl-ST maps correlation coefficient. In fact, just as the \textbf{CC} error does not “see” artificial LOBs for Defl-ST maps, the \textbf{SC} error is not very sensitive to the punctual absence of LOBs, because threshold maps are smooth and their gradients are generally accurate outside LOBs. Otherwise, we can also notice that threshold maps have greater 
\textbf{L2} errors than the Defl-ST maps, illustrating the qualitative observation that threshold maps have a lower overall amplitude than the reference, even though their general shape is respected.

In contrast, the maximal spatio-temporal deflection method on $U$ or on $V$ is able to reproduce real LOB, but also generates artificial LOBs that are impossible to differentiate from real ones. Still, visually, and in terms of correlation coefficient and $L^2$ relative error, the Defl-ST on $U$ method seems to produce the more accurate maps. This result is surprising because it differs from the conclusion of Schuler et al. in \cite{schuler2021reducing}, where the Defl-ST method on $V$ was the only one avoiding artificial LOBs while retaining the highest similarity to the reference.

\section{Discussion and conclusion}
Combined with the Tikonov regularization which has a known smoothing effect, the threshold method tends to generate smooth maps, even in the presence of real discontinuities. Consequently, it would be interesting to combine this method with a more (non-linear) physiological regularization such as the $L^1$ regularization, leading to more accurate reconstructions of the cardiac electrical signals. 

Contrarily to the conclusions of \cite{schuler2021reducing}, our examples show more accurate maps obtained using the Defl-ST on $U$ method rather than the Defl-ST on $V$. This discrepancy may be attributed to the differences in the models and methods of reconstruction of the transmembrane voltage and extracellular potential. For instance, for computing $V$, the authors of \cite{schuler2021reducing} employed the equivalent dipole layer on the epi-endocardial surface while we used a model from \cite{nous} on the pericardium. Thus, the constraint and regularization applied on $V$ was different. Numerous other parameters may also influence the final AT maps, such as the meshes, or the type of temporal smoothing. Henceforth, it is very difficult to recommend one universal method for computing activation maps.



\section*{Acknowledgments}  
This study received financial support from the French Government as part of the “Investments of the Future” program managed by the National Research Agency (ANR), Grant reference ANR-10-IAHU-04.

\bibliography{refs}

\begin{correspondence}
Emma Lagracie\\
200 avenue de la vieille tour, 33400 Talence, France\\
emma.lagracie@inria.fr
\end{correspondence}

\end{document}